 \UseRawInputEncoding

\documentclass[12pt]{amsart}
\usepackage{latexsym, amssymb, amsmath, amscd}
 \usepackage{eucal}

    \newtheorem{rema}{Remark}[section]
    \newtheorem{propo}[rema]{Proposition}

   \newtheorem{theo}[rema]{Theorem}
   \newtheorem{def-theo}[rema]{Definition-Theorem}
 
   \newtheorem{defi}[rema]{Definition}
    \newtheorem{lemma}[rema]{Lemma}
    \newtheorem{corol}[rema]{Corollary}
     
  \newtheorem{rmk}[rema]{Remark}

    \newtheorem{prob}[rema]{Open Problem}

	\newcommand{\nno}{\nonumber}

\newcommand{\olsi}[1]{\,\overline{\!{#1}}} 

\newcommand{\Ker}{\operatorname{Ker\,}}

 \newcommand{\pf}{{\it Proof:}\hspace{2ex}}
 
 \newcommand{\epfv}{\hspace{1em}$\Box$\vspace{1em}}

\newcommand{\bC}{{\mathbb C}}
\newcommand{\bZ}{{\mathbb Z}}

\newcommand{\bQ}{{\mathbb Q}}
\newcommand{\bN}{{\mathbb N}}

\newcommand{\cA}{{\mathcal A}}

\newcommand{\tr}{\operatorname{tr\,}}

\newcommand{\cC}{{\mathcal C}}

\title[Existence of Nonzero Trace-Zero Idempotents]
{Existence of Nonzero Trace-Zero Idempotents in the Group Algebras of Finite Groups}

  \author{Wenhua Zhao and Dan Yan}      
\date{August 09, 2022}
\address{W. Zhao, Department of Mathematics, Illinois State University, Normal, IL 61761. {\it Email}: wzhao@ilstu.edu}
\address{D. Yan, MOE-LCSM, School of Mathematics and Statistics,  Hunan Normal University, Changsha 410081, China. \emph{Email:} yan-dan-hi@163.com }

\begin{document}
\begin{abstract}
Let $G$ be a finite group and $K$ a splitting field of $G$ of characteristic $p>0$. 
Denote by $KG$ the group algebra of $G$ over $K$ 
and $Z(KG)$ the center of $KG$. 
Let $V_G$ be the $K$-subspace of trace-zero elements of $KG$. 
We give some numerical sufficient and necessary conditions for $V_G$ and 
$V_G\cap Z(KG)$, respectively, to be Mathieu subspaces of $KG$ 
in terms of the degrees of irreducible representations 
of $G$ over $K$. The same numerical conditions also characterize 
the finite groups $G$ 
that $KG$ has no nonzero trace-zero idempotents and the finite groups $G$ 
that $KG$ has no nonzero central trace-zero idempotents, respectively. 
\end{abstract}

   




 
\thanks{The first author has been partially supported by the Simons Foundation grant \#278638. The second author is supported by the NSF of China (Grant No. 11871241; 11601146), the China Scholarship Council and the Construct Program of the Key Discipline in Hunan Province.}

 \bibliographystyle{alpha}
    \maketitle


\renewcommand{\theequation}{\thesection.\arabic{equation}}
\renewcommand{\therema}{\thesection.\arabic{rema}}
\setcounter{equation}{0}
\setcounter{rema}{0}
\setcounter{section}{0}

\section{\bf Motivations and the Main Results} \label{S1}

Let $G$ be a finite group, $K$ a splitting field of $G$ 
(unless stated otherwise) and 
$KG$ (or $K[G]$) the group algebra of $G$ over $K$. Denote by $1_G$, or simply $1$, the identity element of $G$ and $KG$. 
We denote by $Z(KG)$ and $V_{KG}$, or simply by $V_G$ if the base field $K$ is clear in context,   
the {\it center} of $KG$ and the $K$-subspace of $KG$ spanned over $K$ by all $1\ne g\in G$, respectively. 
An element $e\in KG$ is said to be an {\it idempotent} 
if $e^2=e$, and a {\it non-trivial idempotent} if $e$ is an idempotent and $e\ne 1$ or $0$. 

For each $u\in KG$, denote by $\operatorname{tr}_G u$, or $\tr u$ if $KG$ is clear in context, the {\it trace} of $u$ 
(i.e., the coefficient of $1_G$ in $u$). 
An element $v\in KG$ is said to be {\it trace-zero} 
if $\tr u=0$ and {\it trace-one} if $\tr u=1$. Hence, the subspace $V_G$ defined above is just 
the subspace of trace-zero elements of $KG$.  

Since an element $e$ is a trace-zero idempotent of $KG$, 
if and only if $1-e$ is a trace-one idempotent of $KG$, 
all the results derived in this paper on trace-zero idempotents of $KG$ 
can be transformed in the obvious ways 
to the results on trace-one idempotents of $KG$.
  
The purpose of this paper is twofold. The first one is to find some  numerical sufficient and necessary conditions in terms of the degrees of irreducible     representations of $G$ over $K$ for that $KG$ has no nonzero trace-zero idempotents and that $KG$ has no nonzero central trace-zero idempotents, respectively. 

To explain the second purpose of the paper, we need to recall the following concept that was introduced by the first author in \cite{GIC} and \cite{MS}.

\begin{defi}
Let $R$ be a commutative ring and $\cA$ a unital $R$-algebra. An $R$-subspace $V$ of $\cA$ is said to be a {\it Mathieu subspace} (MS) of $\cA$ if for all $a, b, c\in \cA$ with $b^m\in V$ for all $m\ge 1$, there exists $N\in \bN$ (depending $a, b$ and $c$) such that $ab^mc\in V$ for all $m\ge N$.  
\end{defi}

A Mathieu subspace  is also called a Mathieu-Zhao space or Mathieu-Zhao subspace 
in the literature (e.g., see \cite{EKC, DEZ}). The introduction of Mathieu subspaces was directly motivated by the Mathieu Conjecture \cite{Ma} and the Image Conjecture \cite{IC}, which both imply the famous Jacobian conjecture \cite{K, BCW, E}. Note that all ideals are Mathieu subspaces, but not conversely. Therefore, the concept of Mathieu subspaces can be viewed as a natural generalization of the concept  of ideals. 

One of the most fundamental theorems in the theory of Mathieu subspaces is the Duistermaat and van der Kallen theorem \cite{DK} which was proved in $1998$.  It solves the abelian case of the Mathieu Conjecture. The theorem in terms of group algebras and Mathieu subspaces  
can be essentially restated as follows.  

\begin{theo}
Let $G=\bZ^n$ (the free abelian group of rank $n$). Then $V_G$ is a MS of 
$\mathbb{C} G$. In terms of $K$-algebra of Laurent polynomials $\bC[x^{-1}, x]$ in $n$ variables, the subspace $V$ of all Laurent polynomials in $x$ with no constant term is a MS of $\bC[x^{-1}, x]$.   
\end{theo}

In order to understand the Duistermaat and van der Kallen theorem in a more general setting, the following problem was proposed 
in \cite[Problem $1.4$]{ZW}.

\begin{prob}
Decide commutative rings $R$ and groups $G$ such that $V_{RG}$ is a Mathieu subspace of the group algebra $RG$ of $G$ over $R$.
\end{prob}

The following proposition and theorem among some other results 
were also proved \cite{ZW}.  

\begin{propo}\label{ZW-Prop}
Let $G$ be a finite group and $K$ an arbitrary field. Then $V_G$ is a MS of $KG$, if and only if $V_G$ does not contain any nonzero idempotent. 
\end{propo}

\begin{theo}\label{ZW-Thm}
Let $G$ be a finite group and $K$ an arbitrary field. Then the following statements hold: 
\begin{enumerate}
  \item[$i)$] If char.\,$K=0$ or char.\,$K=p>|G|$, then $V_G$ 
is a MS of $KG$, where $|G|$ denotes the order of $G$.
  \item[$ii)$] Assume that $G$ is abelian and 
  $|G|=p^rd$ with $r, d\in \bN$ and $p\nmid d$. 
If char.\,$K=p>0$ and $K$ has a $d$-th primitive root of unity, then $V_G$ is a MS of $KG$, if and only if $p>d$. 
\end{enumerate}
\end{theo}

The second purpose of the paper is to find some numerical sufficient and necessary conditions in terms of the degrees of irreducible representations of $G$ over $K$ for that $V_G$ and $V_G\cap Z(KG)$ are Mathieu subspaces 
of $KG$, respectively. 

The two purposes of the paper actually coincide. More precisely, by Proposition \ref{ZW-Prop} 
$V_G$ is a MS of $KG$, if and only if 
$KG$ has no nonzero trace-zero idempotent, and by Lemma \ref{Lma4.1} 
in Section \ref{S4}, $V_G\cap Z(KG)$ is a MS of $KG$, if and only if $KG$ has no nonzero central trace-zero idempotent. 

To state the two main results of the paper, 
we denote by $\phi_i$ $(1\le i\le s)$  the (non-isomorphic) irreducible representations of $G$ over $K$, and $n_i$ the {\it degree} of the representation $\phi_i$. Then the first main theorem of this paper is as follows.

\begin{theo}\label{MainThm}
Let $G$ be a finite group and $K$ a splitting field of $G$ of 
char.\,$K=p>0$. Then $V_G$ is a MS of $KG$, if and only if the following property holds: 
\begin{align} 
c_1 n_1 +c_2n_2+\cdots+c_sn_s \not\equiv 0 \quad (\operatorname{mod}\,\, p)   \tag{$\ast$}
\end{align}
for all integers $0\le c_i\le n_i$ $(1\le i\le s)$ that are 
not all zero. 
\end{theo}

\begin{rmk}\label{ZeroSumLma}
Let $\olsi B$ be a sequence of elements in $\bZ_p$  
formed by $n_i$-copies of $n_i$\,$($mod $p)$ 
with $1\le i\le s$.  
Then it is easy to see that 
the property $(\ast)$ in Theorem \ref{MainThm} for integers 
$n_i$ $(1\le i\le s)$ is equivalent to that the sequence $\olsi B$ is a so-called {\it zero-sum free} sequence $($e.g., see \cite{G}$)$ of  
elements of $\bZ_p$. 
\end{rmk}

%

\begin{rmk}
Given a finite group $G$ and $K$ an arbitrary field. If char.\,$K=0$ or char.\,$K=p>|G|$, then by Proposition \ref{ZW-Prop} and Theorems \ref{ZW-Thm} $KG$ contains no nonzero trace-zero idempotent. If char.\,$K=p\le |G|$, 
let $L$ be a splitting field of $G$, which contains $K$ (e.g., let $L=K[\xi]$, where $\xi$ is a primitive $d$-th root of unity and $d$ is the exponent of $G$. By \cite[Theorem 2.7B, p.\,$58$]{PD} $L$ is a splitting field of $G$). Then we may apply Theorem \ref{MainThm} to decide if $V_{LG}$ is a MS of $LG$. If it is, then by Proposition \ref{ZW-Prop} $V_{KG}$ is a MS of $KG$ and $KG$ contains no nonzero trace-zero idempotent.  
\end{rmk}

In order to state our second main result, we need some basic notations and facts in modular representation theory. Here we freely adapt those fixed and developed in \cite{PD}. 

Assume that char.\,$K=p>0$ and $|G|=p^rd$ with $r, d\in\bN$ and $p\nmid d$.   
Let $\widehat K$ be a field extension of finite degree of the $p$-adic number field $\bQ_p$ such that $X^d-1$ splits completely in $\widehat K[X]$. Then by Lemma \ref{RootLma} proved in Section \ref{S2}, $\widehat K$ is a spitting field for all subgroups of $G$. Let $A$ be the integral closure of the $p$-adic integers in $\widehat K$, and $\pi A$ the unique maximal ideal of $A$ (which contains $pA$). Let $K=A/\pi A$. $K$ is a finite field of characteristic $p>0$. Then by discussion in \cite[Section $2.7$]{PD} $\widehat K$ and $K$ are splitting fields for  all subgroups and quotient groups of $G$. In particular, by \cite[Theorem $1.7B$, p.\,$25$]{PD} all irreducible $\widehat KG$-modules and irreducible $KG$-modules are absolutely irreducible.

Let $f_i$ $(1\le i\le t)$ be the block idempotents of $KG$ and $n_{i, j}$ $(1\le j\le s_i)$ the degrees of all distinct irreducible representations of $G$ over $\widehat K$ which lies in the $i$-th $p$-block of $KG$. 

Our second main theorem of the paper is as follows. 

\begin{theo} \label{MainThm2}
Let $G$ be a finite group and $K$ a splitting field of $G$. Then the following statements hold:
\begin{enumerate}
 \item[$i)$] For every field $L$ with char.\,$L=0$, or  char.\,$L=p>|G|$, $V_{LG}\cap Z(LG)$ is a MS of $LG$.
 
\item[$ii)$] If char.\,$K=p>0$ and $p\nmid |G|$, then 
$V_{KG}\cap Z(KG)$ is a MS of $KG$, if and only if for each 
$\emptyset \ne J\subseteq \{1, 2, ..., s\}$ we have   
\begin{align} \label{MainThm2-eq1}
\sum_{i\in J} n_i^2\not \equiv 0\quad  (\operatorname{mod\,} p).  
\end{align} 

\item[$iii)$] If char.\,$K=p>0$ and $p\,|\, |G|$, then $V_{KG}\cap Z(KG)$ is a MS of $KG$, if and only if  
for all $\emptyset\ne  J\subseteq \{1, 2, ..., t\}$, both the numerator and the denominator of the rational number  
$\frac 1{|G|}\sum_{i\in J} \sum_{j=1}^{s_i} n_{i, j}^2$ in the reduced form are co-prime to $p$.
 \end{enumerate}
\end{theo}

 
Note also that if char.\,$K=0$ or char.\,$K=p>0$ with $p \nmid |G|$, then  by Maschke's Theorem (e.g., see \cite[Theorem $1.3$B]{PD}) $KG$ is semi-simple. For all $1\le i\le t=s$, the dimension over $K$ of $i$-th 
(indecomposable) $p$-block of $KG$ is equal to $n_i^2$.

If char.\,$K=p>0$ with $p \mid |G|$, then for each $1\le i\le t$,  
by \cite[Theorem $4.2$D, p.\,$99$]{PD} and also the fact mentioned above 
on the blocks of $\widehat K G$ (noting char.\,$\widehat K=0)$, 
we see that the dimension of the $i$-th $p$-block of $KG$ over $K$ 
is equal to $\sum_{j=1}^{s_i} n_{i, j}^2$ (in the same notation as in Theorem \ref{MainThm2}). 
Then it is easy to see that the case of 
Theorem \ref{MainThm2} over all splitting base fields 
can be re-stated in the following more uniform way.  

\begin{corol}\label{MainThm2'}
Let $K$ be a splitting field (of arbitrary characteristic) of 
$G$. Let $t$ be the number of (indecomposable) blocks of $KG$
and $d_i$ $(1\le i\le t)$ the dimension over $K$ 
of the $i$-th block of $KG$. 
Then $V_G\cap Z(KG)$ is a MS of $KG$, if and only if  for all 
$\emptyset\ne  J \subseteq \{1, 2, ..., t\}$, both the numerator and the denominator of the rational number  
$\frac1{|G|}\sum_{i\in J} d_i$ in the reduced form are co-prime to $p$. 
%
\end{corol}

Note that by Lemma \ref{Lma4.1} in Section \ref{S4} the sufficient and necessary conditions 
in Theorem \ref{MainThm2} and Corollary \ref{MainThm2'} are also 
sufficient and necessary conditions for the nonexistence 
of nonzero central trace-zero idempotents in $KG$.  \\

{\bf Arrangements}  We mainly give a proof for Theorem \ref{MainThm} 
in Section \ref{S2} and discuss some applications of 
Theorem \ref{MainThm} in Section \ref{S3}. In Section \ref{S4} 
we mainly give a proof for Theorem \ref{MainThm2}. 

\section{\bf Proof of Theorem \ref{MainThm}} \label{S2}

Throughout this section, $G$ stands for a finite group and $K$ a splitting field of $G$ of characteristic $p>0$. We denote by $\phi_i$ $(1\le i\le s)$  all the (non-isomorphic) irreducible representations of $G$ over $K$ and $n_i$ the degree of the representation $\phi_i$ for all $1\le i\le s$. Here we fix the irreducible representation $\phi_i$ $(1\le i\le s)$ as a $K$-algebra homomorphism from  $KG$ to the $n_i \times n_i$ matrix algebra $M_{n_i}(K)$.  Note that by Burnside's Theorem (e.g., see   \cite[Corollary, p.\,$11$]{PD}) each $\phi_i$ $(1\le i\le s)$ is surjective. 

Let $\cA ={\bigoplus}_{i=1}^s M_{n_i}(K)$ and 
$\cA_i=M_{n_i}(K)$ $(1\le i\le s)$. For each $1\le i\le s$ 
we define $\operatorname{tr}_i:\cA\to K$ to the $K$-linear functional whose restriction on the $i$-th component $M_{n_i}(K)$ of $\cA$  is the usual trace function of the $n_i\times n_i$ matrices, 
and on $j$-th component of $\cA$ is equal to $0$ for all $j\ne i$.   

With the notations and settings above we recall the classical Wedderburn Structure Theorem 
(see \cite[Theorem $1.3A$]{PD}) for $KG$. 


\begin{theo}\label{WedderburnThm}
Define $\Phi: KG\to \cA$ be setting $\Phi(\sigma)=\sum_{i=1}^s \phi_i(\sigma)$ for all $\sigma\in KG$. Then  
$\Phi$ is a surjective homomorphism of $K$-algebras and $\Ker \Phi=J(KG)$, where $J(KG)$ stands for the Jacobson radical 
of $KG$.
\end{theo}


Next, we recall the following theorem which clarifies all co-dimensional one MS' of $\cA$.

\begin{theo}$ \label{Kon-Thm} ($\cite[Corollary $4.5$]{Ko}$)$ 
A codimensional one $K$-subspace $V$ of $\cA$ is a MS of $\cA$, 
if and only if $V=\Ker L$ for some $K$-linear functional $L$ of $\cA$ such that $L=\sum_{i=1}^s  c_i\operatorname{tr}_i$, where   
the coefficients $c_i\in K$ $(1\le i\le s)$ 
satisfy the  property: 
{\it for all $\emptyset\ne J\subseteq \{1, 2, ..., s\}$ with  
not all $c_i$ $(i\in J)$ being zero, we have  
$\sum_{i\in J} c_i \ne 0$.}  
\end{theo}

The following proposition classifies all the co-dimensional one MS' $V$ 
of $KG$ with $V\supseteq J(KG)$ (the Jacobson radical of $KG$).

\begin{propo}\label{CoDimOne}
Let $V$ be a co-dimensional one $K$-subspace of $KG$ with $V\supseteq J(KG)$, and $\chi_i$   the character of $G$ obtained from the irreducible representation $\phi_i$ $(1\le i\le s)$.  
Then  $V$ is a MS of $KG$, 
if and only if $V=\Ker \ell$, where $\ell=\sum_{i=1}^s c_i\chi_i$ for some $c_i\in K$ $(1\le i\le s)$ which satisfy the  property: 
{\it for all $\emptyset\ne J\subseteq \{1, 2, ..., s\}$ with  
not all $c_i$ $(i\in J)$ being zero, we have  
$\sum_{i\in J} c_i \ne 0$.}  
\end{propo}
\pf 
Since $\Phi$ by Theorem \ref{WedderburnThm} is surjective and $V\supseteq J(KG)$, there exists a co-dimensional one $K$-subspace $W$ of $\cA$ such that $V=\Phi^{-1}(W)$ (Actually, $W=\Phi(V)$).  
By Theorem \ref{Kon-Thm} (or \cite[Theorem $5.8.1$]{EKC}), 
$W$ is a MS of $\cA$ iff $W=\Ker L$ for some $K$-linear functional $L$ of $\cA$ such that $L=\sum_{i=1}^s c_i\operatorname{tr}_i$ for some $c_i\in K$ 
$(1\le i\le s)$ which satisfy the property in the proposition. 

Let $\Phi^*$ be the dual map of $\Phi$, i.e., $\Phi^*(\psi)=\psi\circ \Phi$ for all $K$-linear functional $\psi$ of $\cA$.
Then by definition of characters we have $\Phi^*(\operatorname{tr}_i)=\chi_i$ for all $1\le  i\le s$. Since $V=\Phi^{-1}(W)$, $W=\Phi(V)$ and $W=\Ker L$, we have $V=\Ker \ell$ with $\ell=\phi^*(L)=\sum_{i=1}^s c_i\chi_i$. 
Since by Theorem \ref{WedderburnThm} and \cite[Proposition $2.7$]{MS}, $V$ is a MS of $KG$, if and only if $W$ is a MS of $\cA$, the proposition follows.  
\epfv

Now, we can prove Theorem \ref{MainThm} for the semisimple case. 

\begin{lemma}\label{SemiSimpleCase}
Assume that $p \nmid |G|$. Then 
\begin{enumerate}
  \item[$i)$] $\Phi(V_G)=\Ker  F$, where $F=\sum_{i=1}^s n_i 
   \operatorname{tr}_i$.
  \item[$ii)$] $V_G=\Ker \ell$, where $\ell=\Phi^*(F)= \sum_{i=1}^s n_i \chi_i$.
  \item[$iii)$] $V_G$ is a MS of $KG$, if and only if $\sum_{i=1}^s c_i n_i  \ne 0$ for all 
 $\vec{0}\ne (c_1, ..., c_s)\in \{ 0, 1, ..., n_1 \}\times \cdots \times \{ 0, 1, ..., n_s \}$. In other words, 
Theorem \ref{MainThm} holds in this case.  
\end{enumerate}  
\end{lemma}
\pf Note first that, since $p\nmid |G|$, by Maschke's Theorem (e.g., see \cite[Theorem $1.3B$]{PD}) 
the group algebra $KG$ is semisimple, i.e., $J(KG)=0$. Then $\Phi$ in Theorem \ref{WedderburnThm} is a $K$-algebra isomorphism.    

$i)$ Let $1_G\ne g \in G$ and $u=\Phi(g)$. Write $u=\sum_{i=1}^s B_i$ with $B_i\in \cA_i\!:=M_{n_i}(K)$ for all $1\le i\le s$. To show $u \in \Ker  F$, we consider the trace of 
the regular action of $g$ on $KG$ and the trace of the regular action of $u$ on $\cA$.  
Since $\Phi$ in this case is an isomorphism of $K$-algebras, we know that these two traces are equal to one another.

First, by choosing elements of $G$ as a basis of $KG$, it is easy to see that the trace of the regular action of $g$ 
on $KG$ is equal to $0$. 

Second, for each $1\le i\le s$, since the matrix $B_i$ annihilates all elements of $\cA_j$ for all 
$j\ne i$, the trace of the regular action of $B_i$ on $\cA$ is equal to the trace of the regular action of $B_i$ on 
$\cA_i=M_{n_i}(K)$, which can be readily checked to be 
$n_i \operatorname{tr}_i (B_i)$ by choosing the matrices $E_{k, \ell}$'s  
(with $1$ as the $(k, \ell)$-entry and zero elsewhere) as a basis of $\mathcal A_i=M_{n_i}(K)$. Hence the trace of 
the regular action of $u$ on $\cA$ is equal to 
$\sum_{i=1}^s n_i \operatorname{tr}_i(B_i)$. 
Therefore $\sum_{i=1}^s n_i \operatorname{tr}_i(B_i)=0$, i.e., 
$u\in \Ker F$, whence $\Phi(V_G)\subseteq \Ker F$ for $V_G$ is spanned over $K$ by all non-identity elements of 
$G$. Since $\Phi(V_G)$ and $\Ker F$ are both co-dimensional 
one subspaces of $\cA$, we have $\Phi(V_G)= \Ker F$, and statement $i)$ follows. 

$ii)$ follows immediately from $i)$, and the fact that $\Phi^*(\operatorname{tr}_i)=\chi_i$ for all $1\le i\le s$. 
$iii)$ follows immediately from $ii)$, Proposition \ref{CoDimOne} above and the fact that $n_i\ne 0$ for all $1\le i\le s$. \epfv  

Next we shall give a proof of Theorem \ref{MainThm} for the case $p\,|\, |G|$. In order to do so we need first to prove some lemmas, for which we do not assume $p\,|\, |G|$.

\begin{lemma}\label{OrthogonalIdem}
Assume that $KG$ has $k$ nonzero (pairwise) orthogonal idempotents $e_i$ $(1\le i\le k)$ with $k\ge p$. 
Then $V_G$ is not a MS of $KG$.   
\end{lemma}

\pf For each $1\le j\le k$ set $s_j=\sum_{i=1}^j e_i$. Then it is easy to see that $s_j$ $(1\le j\le k)$ are all nonzero idempotents. By \cite[Theorem $3.5$ in p.\,$48$]{P}, 
which was first proved by Zalesskii \cite{Za} in $1972$, 
$\tr s_j\in \bZ_p=\bZ/p\bZ$ for all $1\le j\le k$. 
Since $k\ge p$, there exists either $1\le j_0\le k$ such that 
$\tr s_{j_0}=0$, or $a<b$ such that 
$\tr s_a=\tr  s_b$, whence 
$\tr(s_b-s_a)=0$. Since $s_b-s_a=e_{a+1}+\cdots+e_b$ is also a nonzero idempotent, $V_G$ in either case contains 
a nonzero trace-zero idempotent, namely, $s_{j_0}$ or $s_b-s_a$.
Then by Proposition \ref{ZW-Prop}  
$V_G$ is not a MS of $KG$. 
\epfv 
 
\begin{lemma}\label{NormalLma}
Assume that $V_G$ is a MS  of $KG$. Then 
\begin{enumerate}
\item[$i)$] $p>\sum_{i=1}^s n_i$.
\item[$ii)$]  $G$ has a normal Sylow $p$-subgroup.
\end{enumerate}  
\end{lemma}
\pf $i)$ Let $m=\sum_{i=1}^s n_i$ and assume $p\le m$. 
We view $\cA$ as a $K$-subalgebra of the matrix algebra $M_m(K)$ by putting $M_{n_i}(K)$ $(1\le i\le s)$ as the $(i, i)$-th block on the diagonal. Let $f_i$ $(1\le i\le m)$ be the $m\times m$ diagonal  matrix with all entries being zero except the $(i, i)$-th entry being $1$. Hence $f_i$ $(1\le i\le m)$ are orthogonal idempotents of $\cA$ such that $\sum_{i=1}^m f_i=1_\cA$. 
Note that by Theorem \ref{WedderburnThm}, $\Ker \phi=J(KG)$ and hence is nilpotent. 
Then by \cite[Lemma $3.7$\,(i), p.\,$49$]{P} there exist nonzero orthogonal idempotents $e_i \in KG$ $(1\le i\le m)$ such that $\Phi(e_i)=f_i$. Hence, 
by Lemma \ref{OrthogonalIdem}, $V_G$ is not a MS  of $KG$. Contradiction. 

$ii)$ By $i)$ we have $p>n_i$ for each $1\le i\le s$, 
whence $p\nmid n_i$. It is well known that 
$n_i$ $(1\le i\le s)$ are also the degrees of 
irreducible Brauer characters of $G$ (by the very definition of the irreducible Brauer characters, e.g., see \cite[p.\,$17$]{N}). 
Then by \cite[Theorem $2.4$]{M2}
\footnote{The proof of \cite[Theorem $2.4$]{M2} used the classification of finite simple groups. For a proof of this theorem without using the classification, see \cite[Theorem on p.\,$121$]{Manz}. For a short proof for a special case of the ``only if" part of the theorem, see \cite{R}. } 
(which is a modular version of the It\^o-Michler Theorem \cite{I, M2}, see also \cite{N2}) 
%
$G$ has a normal Sylow $p$-subgroup. 
\epfv

Recall that an element $g\in G$ is said to be a {\it $p$-singular} element, or simply a {\it $p$-element} of $G$ 
if the order of $g$ is a power of $p$. $g$ is said to be a {\it $p$-regular element}, or simply a {\it $p'$-element} of $G$ of the order of $g$ is  co-prime to $p$. The identity element $1_G$ is the only both $p$-singular and $p$-regular element of $G$. 
  
It is easy to check that every element $g\in G$ can be written as the product of a $p$-element and $p'$-element, 
which commute with one another. We denote by $G^\circ$ the set of all $p'$-elements of $G$.

\begin{lemma} \label{SumZeroLma}
Let $e$ be an idempotent of $KG$ and $h\in G\backslash G^\circ$, i.e., $h$ is not a $p'$-element. 
Then the sum of the coefficients in $e$ of the elements in the conjugacy class $\cC$ of $h$ is equal to $0$.
\end{lemma}
\pf  Write $|G| = p^r d$ for some $r, d\in \bN$ with $p \nmid d$. Then it is easy to see that $g^{p^r}\in G^\circ$ for all $g\in G$. Write the idempotent $e$ as  
\begin{eqnarray}
e=\sum_{a\in G\backslash G^\circ}c_a \,a +\sum_{b \in G^\circ}c_b \, b \label{eq2.1}
\end{eqnarray}
with $c_a, c_b\in K$ $(a\in G\backslash G^\circ$, $b \in G^\circ)$.  
 
By \cite[Lemma $3.1$, p.\,$45$]{P} there exists $\beta$ in the commutator subspace $[KG, KG]$ of $KG$ such that   
\begin{align}
e=e^{p^r}=\sum_{a\in G\backslash G^\circ} c_a^{p^r} a^{p^r}+\sum_{b \in G^\circ}c_b^{p^r}b^{p^r}+\beta. \nno
\end{align}
Since $a^{p^r}\in G^\circ$ for all $a\in  G\backslash G^\circ$, we further have 
\begin{align}
 e=e^{p^r}
 = \sum_{b \in G^\circ}d_b b +\beta  \label{eq2.2}
\end{align}
for some $d_b\in K$ $(b\in G^\circ)$.

Then by eqs.\,(\ref{eq2.1}) and (\ref{eq2.2}) we get  
\begin{align}
 \sum_{a\in G\backslash G^\circ}c_a \,a +\sum_{b \in G^\circ}(c_b-d_b) \, b =\beta \in [KG, KG]. \label{eq2.3}
\end{align}
Since $h \in G\backslash G^\circ$, its conjugacy class $\cC \subset G\backslash G^\circ$. Then by \cite[Lemma 3.2, p.\,$46$]{P} and the equation above we have 
$\left(\sum_{a\in\mathcal C} c_a\right)^{p^r}=\sum_{a\in\mathcal C} c_a^{p^r}=0$, whence the lemma follows. 
\epfv


The next lemma will play a crucial role in the proof of Theorem \ref{MainThm}.

\begin{lemma}\label{G=G/H}
Let  $H$ be a normal $p$-subgroup of $G$. Then $V_G$ is a MS  of $KG$, if and only if $V_{G/H}$ is a MS of the group algebra 
$K[G/H]$ of $G/H$ over $K$.
\end{lemma}
\pf  
Let $\rho: K[G]\to K[G/H]$ be the $K$-algebra homomorphism induced by the quotient map from $G$ to $G/H$. 
Then by \cite[Lemma 1.8, p.\,$10$]{P} $\Ker \rho =\omega (KH)KG$, where $\omega (KH)$ is the augmentation ideal of $KH$ which is spanned by $(h-1)$ $(h\in H)$ over $K$. Therefore $K[G/H]\simeq KG/\omega (KH)KG$. 

Let $g_i$ $(1\le i\le k\!:=|G/H|)$ with $g_1=1_G$ be the representatives of the distinct (right) cosets of $H$ in $G$. Then for every idempotent $e\in KG$, there exist some $a_i \in KH$ $(1\le i\le k)$ such that 
\begin{align}\label{G=G/H-peq1}
e=a_1 1_G+a_2 g_2+\cdots+a_k g_k.
\end{align}
Hence, the coefficient of $1_G$ in $e$ is equal to the coefficient of $1_H$ in $a_1$, i.e., $\operatorname{tr}_G e=\operatorname{tr}_H a_1$. Furthermore, since $H\unlhd G$, the $G$-conjugacy class of each $h\in H$ is contained in $H$. 
Since $H$ is a $p$-group, by Lemma \ref{SumZeroLma} the sum of the coefficients in $a_1$ over the $G$-conjugacy 
class of each $1_H\ne h\in H$ is equal to zero, 
whence $\rho(a_1)=\operatorname{tr}_H a_1$. 

On the other hand, by eq.\,(\ref{G=G/H-peq1}) and 
the definition of $\rho$ we also have 
$\operatorname{tr}_{G/H} \rho(e)= \rho(a_1)$. 
In summary, the following equations hold:
\begin{align}\label{G=G/H-peq2}		
\operatorname{tr}_G e=\operatorname{tr}_H a_1=\rho(a_1)=\operatorname{tr}_{G/H} \rho(e).
\end{align}
    
$(\Rightarrow)$ Assume that $V_G$ is a MS  of $KG$ and $V_{G/H}$ 
is not a MS of $K[G/H]$. Then by Proposition \ref{ZW-Prop} there 
exists a nonzero idempotent $u\in V_{G/H}$ (with $\operatorname{tr}_{G/H} u=0$).
Since $H\unlhd G$, it follows easily from  
\cite[Lemma 1.6, p.\,$70$]{P} that 
$\Ker \rho$ ($=\omega (KH)KG$) is a nilpotent ideal of $KG$. Then   
by \cite[Lemma 3.7 (i), p.\,$49$]{P} there exists an idempotent $e\in KG$ 
such that $\rho (e)=u$. 
By eq.\,(\ref{G=G/H-peq2}) $\operatorname{tr}_G e=\operatorname{tr}_{G/H} u=0$, i.e, $e\in V_G$. Since $V_G$ is MS of $KG$, by Proposition \ref{ZW-Prop} we have $e=0$, whence $u=0$. Contradiction. 
 
$(\Leftarrow)$ Assume that $V_{G/H}$ is a MS of $K[G/H]$ and 
$V_G$ is not a MS  of $KG$. Then by Proposition \ref{ZW-Prop} there 
exists a nonzero idempotent $e\in V_G$ (with $\operatorname{tr}_G e=0$).
Note that $\rho(e)\neq 0$, for $\Ker \rho=\omega(KH)KG$ is nilpotent (as pointed out above). Therefore $\rho(e)$ is a nonzero idempotent 
of $K[G/H]$, and by eq.\,(\ref{G=G/H-peq2}) $\operatorname{tr}_{G/H} \rho(e)=0$, whence $\rho(e)\in V_{G/H}$.  
Since $V_{G/H}$ is a MS of $K[G/H]$, 
by Proposition \ref{ZW-Prop} $\rho(e)=0$. 
Contradiction. \quad
$\square$ 
 
\begin{rmk} We will show later in Proposition \ref{Closed-H-G/H} that the 
$p$-subgroup condition on $H$ for the ``only if" part of Lemma \ref{G=G/H} can actually be dropped.
\end{rmk}

Now we can prove Theorem \ref{MainThm} as follows.\\

{\bf Proof of Theorem \ref{MainThm}}: 
$(\Rightarrow)$ Assume that $V_G$ is a MS  of $KG$. Then by Lemma \ref{NormalLma} $ii)$ $G$ has a normal Sylow $p$-subgroup $H$. By Lemma \ref{G=G/H} $V_{G/H}$ is a MS  of $K[G/H]$. Since $p\nmid|G/H|$, it follows from Lemma \ref{SemiSimpleCase} $iii)$ that the property $(\ast)$ in Theorem \ref{MainThm} holds for $K[G/H]$. 

On the other hand, by \cite[Proposition $4.7$, p.\,67]{P} the Jacobson radical $J(KG)$ of $KG$ is equal to the ideal 
$I_H$ of $KG$ generated by $(h-1)$ $(h\in H)$.  
Since by Theorem \ref{WedderburnThm}  
$K[G/H]\simeq KG/I_H= KG/J(KG)$,  
the degrees of irreducible representations 
of $KG$ over $K$ (or equivalently, irreducible representations of $G$ over $K$) are the same as the degrees of irreducible representations 
of $K[G/H]$ over $K$ (or equivalently, irreducible representations of $G/H$ over $K$). Therefore, the property $(\ast)$ in Theorem \ref{MainThm} holds also for $KG$.

$(\Leftarrow)$ Assume that the property $(\ast)$ in Theorem \ref{MainThm} holds for $KG$. If $p \le n_{i_0}$ for some $1\le i_0\le s$, then the property $(\ast)$ fails for the coefficients 
$c_{i_0}=p$ and $c_i=0$ $(i\ne i_0)$. Therefore,    
$p>n_i$, and hence $p\nmid n_i$, for all $1\le i\le s$. 
Then by \cite[Theorem $2.4$]{M2} or \cite[Theorem on p.\,$121$]{Manz} or the theorem in the short note \cite{R}, 
%
$G$ has a normal Sylow $p$-subgroup $H$. Since  
as pointed above, the degrees $n_i$ $(1\le i\le s)$ of irreducible representations of $G$ over $K$ are the same as the degrees of irreducible representations of $G/H$ over $K$, the property $(\ast)$ in Theorem \ref{MainThm} 
holds for $K[G/H]$. Since $p\nmid |G/H|$, by Lemma \ref{SemiSimpleCase} $iii)$ $V_{G/H}$ is a MS of $K[G/H]$, and by Lemma \ref{G=G/H} 
$V_G$ is a MS of $KG$.  
\epfv 
 
To end this section we point out that Theorem \ref{ZW-Thm}, 
$ii)$ also follows from Theorem \ref{MainThm}. 
First, by the lemma given below, the condition 
$p>d$ implies that $K$ is a splitting field of 
$G$. Second, it is well known  
(e.g., see \cite[Corollary on p.\,$19$]{PD}) that all irreducible representations of an abelian group are of dimension one, 
i.e., $n_i=1$ for all $1\le i\le s$.  
Then it is easy to see that the property $(\ast)$ 
in Theorem \ref{MainThm} holds, if and only if $p>s$.

On the other hand, by \cite[Theorem $1.5$]{PD} 
the number of irreducible representations of $G$ over $K$ 
is equal to the number of conjugacy classes 
of $p'$-elements of $G$, i.e., $s=d$ (since $G$ is abelian). 
Then by the facts above and also the lemma below, 
Theorem \ref{ZW-Thm} $ii)$ follows from Theorem \ref{MainThm}.

\begin{lemma}\label{RootLma}
Assume that $|G|=p^rd$ for some $r, d\in \bN$ such that $p\nmid d$. 
Let $L$ be a field of characteristic $p>0$, which contains a primitive $d$-th root of unity. 
Then $L$ is a splitting field of $G$.  
\end{lemma}
\pf Since $L$ contains a primitive $d$-th root of unity, the polynomial $X^d-1$ splits into the product of some degree one polynomials in $L[X]$. Let $n=|G|$. Since $X^n-1=X^{p^rd}-1=(X^d-1)^{p^r}$, the polynomial $X^n-1$ also splits into the product of degree one polynomials in $L[X]$, hence so does $X^m-1$, where $m$ is the exponent of 
$G$. Then by \cite[Theorem 2.7B, p.\,$58$]{PD} 
$L$ is a splitting field of $G$. \quad $\square$  
 \vspace{3mm}

\section{\bf Some Applications of Theorem \ref{MainThm}}
\label{S3}

Throughout this section $G$ stands for a finite group with $|G|=p^rd$ and $ p \nmid d$, and $K$ a splitting filed of $G$ of characteristic $p>0$. We will also freely use the notation fixed in Sections \ref{S1} and \ref{S2}, e.g., $n_i$ $(1\le i\le s)$ stand for the degrees 
of all irreducible representations of $G$ over $K$.

\begin{propo}\label{p>|G/H|}
Assume that $|G|=p^r d$ with $p\nmid d$.  If $p>\sum_{i=1}^s n_i^2$, 
then $V_G$ is  a MS  of $KG$. In particular, if $G$ has a normal Sylow $p$-subgroup and $p>d$, then $V_G$ is a MS of $KG$. 
\end{propo}
\pf First, it is easy to see that the condition $p>\sum_{i=1}^s n_i^2$ implies the property $(\ast)$ in Theorem \ref{MainThm}, and hence by the same theorem $V_G$ is a MS  of $KG$. 

Now assume that $G$ has a normal Sylow $p$-subgroup $H$. Then,   
as pointed out in the proof of Lemma \ref{G=G/H} in Section \ref{S2}, $n_i$ $(1\le i\le s)$ are also the degrees 
of all irreducible representations of $G/H$ over $K$. 
Since $p\nmid d=|G/H|$, by Maschke's Theorem (e.g., see \cite[Theorem $1.3B, \, ii)$]{PD} $K[G/H]$ is semi-simple. In particular, we have $\sum_{i=1}^s n_i^2=d$ 
by Theorem \ref{WedderburnThm}, and hence $p>\sum_{i=1}^s n_i^2$. Then by the first part of the proposition proved above, $V_G$ is a MS  of $KG$.
\quad $\square$


\begin{rmk} \label{p-grp} 
Let $G$ be a $p$-group, then $V_G$ by the proposition above is a MS of $KG$.
Actually, by \cite[Proposition $4.7$, p.\,67]{P} $KG$ is a 
local $K$-algebra and hence, 
does not have any non-trivial idempotents. 
Then by \cite[Theorem $4.2$]{MS} 
every $K$-subspace $V$ of $KG$ such that $1\not\in V$ 
is a MS of $KG$. 
\end{rmk}

The next two propositions claim that the set of all 
finite groups $G$ over some fixed splitting 
fields $K$ is closed under 
taking subgroups, and taking   
quotient groups as well as taking semi-product with $p$-groups.

\begin{propo}\label{Closed-H-G/H}
Let $Q$ be a subgroup of $G$ or a quotient group of $G$. 
Assume that $K$ is also a splitting field of $Q$ and $V_G$ is a MS of $KG$. 
Then $V_Q$ is a MS of $KQ$. 
 
\end{propo} 
\pf The case that $Q$ is a subgroup can be seen easily from the assumption and Proposition \ref{ZW-Prop}. So we assume that  
$Q=G/H$ for some $H\unlhd G$.

Let $\rho: KG\to K[G/H]$ be the $K$-algebra homomorphism induced by the quotient map from $G$ to $G/H$ and $\psi_i: K[G/H]\to M_{n_i}(K)$ $(1\le  i\le t)$ 
all (non-isomorphic) irreducible representations of 
$G/H$ over $K$. 
Then by Burnside's Theorem (e.g., see   \cite[Corollary, p.\,$11$]{PD}) each $\psi_i$ is surjective. Since $\pi$ is surjective, so is the composition $\phi_i\!:=\psi_i\circ \rho$.  Hence $\phi_i$ $(1\le i\le t)$ are non-isomorphic irreducible representations of $G$ over $K$ 
of degree $n_i$. Therefore, the condition $(\ast)$ in Theorem \ref{MainThm} for the degrees of irreducible modules of $G$ over $K$ 
implies the same condition $(\ast)$ for the degrees 
of irreducible modules of $G/H$ over $K$ 
(by letting $c_i=0$ for all $t+1\le i\le s$). 
Then by applying Theorem \ref{MainThm} to $G/H$ 
the proposition follows.
\epfv
 
A special case of the proposition above for quotient groups of $G$   follows also directly from Proposition \ref{ZW-Prop} and the lemma below which has some interests on its own.

\begin{lemma}\label{p'LiftingLma}
Let $H\unlhd G$ such that $p\nmid |H|$ and $\rho: KG\to K[G/H]$ be the $K$-algebra homomorphism induced by the quotient map from $G$ to $G/H$. 
Then for every idempotent $e \in K[G/H]$ there exists an idempotent 
$u\in KG$ such that $\rho(u)=e$ and  
$\operatorname{tr}_G u=\frac1{|H|}\operatorname{tr}_{G/H} e$.
\end{lemma}  

\pf Let $k=|G/H|$ and $g_i$ $(1\le i\le k)$ be 
the representatives of all the cosets of 
$H$ in $G$ with $g_1=1_G$. 
For each $g\in G$ we set also 
$\bar g=\rho(g)$.  
Write $e=\sum_{i=1}^k c_i \bar g_i$ 
for some $c_i\in K$ $(1\le i\le k)$. 
Let $E_H=\frac1{|H|} \sum_{h\in H} h$ 
(i.e., the principal idempotent of 
the group algebra $KH$ of $H$). 
Then it is well known (and also 
easy to check) that we have
\begin{align}\label{p'LiftingLma-peq1}
 E_H^2=E_H \quad \text{and} \quad E_Hh=E_H
\end{align}
for all $h\in H$. Furthermore, 
by using the assumption $H\unlhd G$ 
we have also that $ E_Hg=gE_H$ for all $g\in G$, whence 
$E_H\in Z(KG)$. Set 
\begin{align}\label{p'LiftingLma-peq2}
v\!:=\sum_{i=1}^k c_i g_i \quad \text{and}\quad u\!:=E_Hv.
\end{align}
Then $\rho(v)=e$, $\rho(E_H)=1$, and hence $\rho(u)=e$. Consequently, $v^2-v\in \Ker\rho$, for $\rho(v^2)=\rho(v)^2=e^2=e=\rho(v)$, 
and $v^2=v+w$ for some $w\in \Ker\rho$. Furthermore, 
by \cite[Lemma $1.8$, p.\,$10$]{P} there exist some 
$\alpha_h \in KG$ $(h\in H)$ such that 
$w=\sum_{h\in H}(h-1)\alpha_h$. Now, by eq.\,(\ref{p'LiftingLma-peq1}) 
and the fact $E_H\in Z(KG)$ we have 
\begin{align*}
u^2&=E_H v E_H v=E_H^2v^2=E_H v^2=E_H(v+w)=E_H(v+\sum_{h\in H}(h-1)\alpha_h)\\
 &= E_Hv+ \sum_{h\in H}(E_h h-h)\alpha_h=E_Hv=u.
\end{align*}
Therefore $u$ is an idempotent of $KG$ with $\rho(u)=e$, and by 
eq.\,(\ref{p'LiftingLma-peq2})  $\operatorname{tr}_G u=\frac1{|H|}\operatorname{tr}_{G/H} e$, as desired. 
\epfv

\begin{propo}\label{Semi-Prod}
Let $P$ be a $p$-group and $\widetilde{G}=P\rtimes G$. 
Then $V_G$ is a MS of $KG$, if and only if  
$V_{\widetilde{G}}$ is a MS 
of $K\widetilde{G}$.
\end{propo}
\pf Since $\tilde{G}/P\cong G$ and $P$ is a normal $p$-subgroup of $\widetilde{G}$, the proposition follows immediately from Lemma \ref{G=G/H}.
\epfv 

\begin{propo}\label{VonNewmann}
Assume that $K$ is splitting field for all subgroups and quotient groups of $G$ 
and $V_G$ is a MS  of $KG$. 
Then for every subgroup or quotient $Q$, the following statements hold: 
\begin{enumerate}
\item[$i)$] the degree of every irreducible module of $Q$ over $K$ is less than $p$.
\item[$ii)$] the number of $p$-blocks of $KQ$ is less than $p$.
\end{enumerate}
\end{propo}
\pf By Proposition  \ref{Closed-H-G/H} we may assume $Q=G$. Then 
$i)$ follows from Lemma \ref{NormalLma} $i)$. Since the number of $p$-blocks of $KG$ is equal to the number of (orthogonal and primitive) 
block idempotents of $KG$, statement $ii)$ follows immediately from Lemma \ref{OrthogonalIdem}.
%
\epfv


Next, we consider the extreme case that char.\,$K=2$.

\begin{propo}\label{2-grps}
Let $L$ be a field of characteristic $p=2$, which is a splitting field of all subgroups of $G$ (e.g., when the polynomial $X^d-1$ splits completely in $L[X]$, see Lemma \ref{RootLma}). Then the following statements are equivalent: 
\begin{enumerate}
\item[$i)$] $V_G$ is a MS  of $LG$. 
\item[$ii)$] $LG$ does not contain any non-trivial idempotent.
\item[$iii)$] $G$ is a $2$-group.
\end{enumerate} 
\end{propo}
\pf Let $e$ be an idempotent of $LG$. Then by   \cite[Theorem $3.5$ in p.\,$48$]{P}, which was first proved by Zalesskii \cite{Za} in $1972$,  
$\tr e=0$, or $1$. Hence the equivalence of statements $i)$ and $ii)$ follows immediately from Proposition \ref{ZW-Prop}. 
The  $iii)\Rightarrow i)$ part follows from Remark \ref{p-grp}. 
So it suffices to show $i)\Rightarrow iii)$.

Assume that $i)$ holds but $iii)$ fails. By Cauchy's theorem there exists $g\in G$ of odd order $k\ge 3$. 
Let $H$ be the cyclic subgroup of $G$ generated by $g$. 
Then by Proposition \ref{Closed-H-G/H}, $V_H$ is a MS  of $LH$, and  
the condition $(\ast)$ in Theorem \ref{MainThm} 
holds for the degrees of irreducible representations of $H$, 
which together with the assumption  
$p=2$ implies that $H$ has only one irreducible representation, 
namely, the trivial one.

On the other hand, since $2\nmid k=|H|$, by \cite[Theorem $1.5$]{PD} the number of irreducible modules of $H$ over $K$ is equal to the number 
of conjugacy classes of odd order elements of $H$, 
which is $k$ since $H$ is abelian. Therefore $k=1$, 
contradiction.
\epfv  
 
Finally, from Theorem \ref{MainThm} and 
Lemma \ref{NormalLma} we see that, for the degrees $n_i$ 
$(1\le i\le s)$ of irreducible representations of $G$ over $K$, the property $(\ast)$ in Theorem \ref{MainThm} implies that $p>\sum_{i=1}^s n_i$. 
Actually, this implication holds for all finite sequences of positive integers, as we will show in the next corollary.  

\begin{corol}
Let $n_i$ $(1\le i\le s)$ be arbitrary positive integers which 
satisfy the property $(\ast)$ in Theorem \ref{MainThm}. Then  
$p>\sum_{i=1}^s n_i$.
\end{corol}
\pf Let $\olsi B$ be a sequence of elements in $\bZ_p$  
formed by $n_i$-copies of $n_i$\,$($mod $p)$ 
with $1\le i\le s$. Then as pointed out in Remark \ref{ZeroSumLma} 
the property $(\ast)$ in Theorem \ref{MainThm} for integers 
$n_i$ $(1\le i\le s)$ is equivalent to that $\olsi B$ is a 
{\it zero-sum free} sequence of elements of $\bZ_p$. 
Since the sequence $\olsi B$ is of length $\sum_{i=1}^s n_i$, 
the corollary follows from \cite[Theorem $10.2$]{G} for 
the abelian group $\bZ_p$. 
\epfv

The corollary above can also be proved directly by applying similar arguments as those given in the proof of Lemma \ref{OrthogonalIdem} to 
the sequence of partial sums of $\olsi B$ in the proof above.

\section{\bf Proof of Theorem \ref{MainThm2}}\label{S4}

Throughout this section $G$ stands for a finite group and $K$ a field 
(of arbitrary characteristic). Let $f_i$ $(1\le i\le t)$ be all the (central) block idempotents of $KG$. The notations fixed in 
the previous sections are still in force, e.g., $|G|=p^rd$ 
with $p\nmid d$.   

In this section we mainly give a proof for Theorem \ref{MainThm2}. 

\begin{lemma}\label{Lma4.1}
The following statements are equivalent: 
\begin{enumerate}
\item[$i)$]   $V_G\cap Z(KG)$ is a MS  of $KG$;
\item[$ii)$]   $KG$ does not contain any nonzero central trace-zero  idempotent.  
\item[$iii)$]  $\sum_{i\in J} \tr f_i \neq 0$ for all non-empty 
 $J\subseteq \{1, 2, ..., t\}$. 
\end{enumerate}
\end{lemma}
\pf $i)\Rightarrow ii)$ Assume that  $KG$ contains a nonzero central idempotent $e$ with $\tr e=0$. Let $g\in G$ such that the coefficient 
$c_g$ of $g$ in $e$ is nonzero. Then $\tr (g^{-1}e)=c_g\ne 0$,  
whence $g^{-1}e\not\in V_G\cap Z(KG)$. Since $KG$ is of finite dimension over $K$ and hence algebraic over $K$, by \cite[Theorem $4.2$]{MS}
$V_G\cap Z(KG)$ is not a MS of $KG$.  Contradiction. 
 
The $ii)\Rightarrow i)$ part follows immediately from \cite[Theorem $4.2$]{MS}.

$ii)\Leftrightarrow iii)$ Since $f_i$ $(1\le i\le t)$ are all primitive idempotents of $Z(KG)$, by 
\cite[Theorem $3.11$, p.\,$55$]{N} every idempotent of $Z(KG)$ has the form $\sum_{i\in J} f_i$ for some non-empty 
$J\subseteq \{1, 2, ..., t\}$. Hence statements $ii)$ and $iii)$ are equivalent to one another.
$\square$

\begin{rmk}
By a similar argument as the proof of \cite[Proposition $3.1$]{ZW}, one can actually show that, replacing the base field $K$ by any commutative base ring $R$, $V_{RG}\cap Z(RG)$ is a MS of $RG$, if and only if every 
$\beta\in V_{RG}\cap Z(RG)$ with $\beta^m\in V_{RG}\cap Z(RG)$ for all $m\ge 1$ is nilpotent.
\end{rmk}

\begin{corol}\label{Corol4.2}
If $V_G$ is a MS of $KG$, then so does $V_G\cap Z(KG)$.
\end{corol}
\pf Since $V_G$ is a MS of $KG$, by Proposition \ref{ZW-Prop} $V_G$ does not contain any nonzero idempotent of $KG$. Hence, neither does $V_G\cap Z(KG)$. Then by Lemma \ref{Lma4.1}, $V_G\cap Z(KG)$ is a MS of $KG$.
\epfv


\begin{corol}
Assume that char.\,$K=p>0$ with $p\mid |G|$, and $V_G\cap Z(KG)$ is a MS  of $KG$. Then every block idempotent $f_i$ $(1\le i\le t)$ is of full defect.
\end{corol}
\pf For each $1\le i\le t$, let $r_i$ be the defect of the block idempotent $f_i$. Then by \cite[Theorem $2.1\, (b)$]{M1} 
$\tr f_i = p^{r-r_i}\frac xy$ (as elements of $K$) 
for some $x, y\in \bN$ which are co-prime to $p$. 
If $r_i<r$ for some $1\le i\le t$, then statement $iii)$ in Lemma \ref{Lma4.1} fails for the subset $J=\{i\}$, and by the same lemma, $V_G\cap Z(KG)$ is not a MS  of $KG$. Contradiction. 
\epfv

Now we can give a proof for Theorem \ref{MainThm2}. We will freely use the notation and conventions introduced in Section \ref{S1} related with  Theorem \ref{MainThm2} and Corollary \ref{MainThm2'}. 
\\
 
{\bf Proof of Theorem \ref{MainThm2}:} $i)$ By Theorem \ref{ZW-Thm} $i)$, $V_G$ is a MS of $LG$, and by Corollary \ref{Corol4.2}, $V_G\cap Z(LG)$ is a MS of $LG$.
 
$ii)$ Assume that char.\,$K=p>0$ with $p\nmid |G|$, 
then by \cite[Theorem $3.5, \, ii)$, p.\,$75$]{PD}, 
$s_i=1$ and $n_{i, 1}=n_i$ for all $1\le i\le t$. 
Consequently, we have $t=s$, and by \cite[eq.\,$(8)$, p.\,93]{PD} and \cite[Theorems $3.5$ and $4.2$D]{PD} $\tr f_i=\frac 1{|G|} n_i^2$ (as an element of $K$) for all $1\le i\le s=t$. Then by Lemma \ref{Lma4.1} $iii)$, 
$V_G\cap Z(KG)$ is a MS of $KG$, if and only if for each non-empty $J\subseteq \{1, 2, ..., s\}$ we have
$$
\frac1{|G|}\sum_{i\in J} n_i^2 \ne 0_K. 
$$
Since $|G|\ne 0_K$ and $n_i^2\in \bN$ for all $1\le i\le s$, we see 
that the inequality above holds, if and only if 
$\sum_{i\in J} n_i^2 \not\equiv  0$ (mod $p$). 
Hence statement $ii)$ follows.

$iii)$ By \cite[Theorem $4.2D$, p.\,99]{PD} (noting that the field $K$ there is our $\widehat K$), $\tr\!f_i=\frac 1{|G|} \left( \sum_{i=1}^{s_i} n_{i, j}^2\right)$ (as an element of $K$) for all $1\le i\le t$. 
Then statement $iii)$ follows immediately from Lemma \ref{Lma4.1} 
$iii)$. \epfv

{\bf Acknowledgments} The authors are grateful to an anonymous referee for some valuable suggestions on an earlier version of the paper. The first author also thanks Sunil Chebolu for some personal communications.


\begin{thebibliography}{FLM2}


\bibitem[BCW]{BCW} H. Bass, E. Connell and D. Wright, {\it The Jacobian Conjecture, Reduction of Degree and Formal Expansion of the Inverse}. Bull.  Amer. Math.  Soc.  \textbf{7}, (1982), 287--330. 

\bibitem[DEZ]{DEZ} H. Derksen, A. van den Essen and W. Zhao, {\it The Gaussian Moments Conjecture and the Jacobian Conjecture}.  Israel J. Math. {\bf 219} (2017), no.\,2, 917--928. See also arXiv:1506.05192 [math.AC]. 

\bibitem[DK]{DK}  J. J. Duistermaat and W. van der Kallen, {\it Constant Terms in Powers of a Laurent Polynomial}. Indag. Math. (N.S.) {\bf 9}\, (1998), no. 2, 221–231.

\bibitem[E]{E} A. van den Essen, {\it Polynomial Automorphisms and the Jacobian Conjecture}. Prog. Math., Vol.190, Birkh\"auser Verlag, Basel, 2000. 

\bibitem[EKC]{EKC} A. van den Essen, S. Kuroda and A. J. Crachiola, {\it Polynomial Automorphisms
and the Jacobian Conjecture: New Results from the Beginning of the 21st
Century}. Frontiers in Mathematics. Birkh\"auser $2021$. 

\bibitem[G]{G} D. J. Grynkiewicz, {\it  Structural additive theory}.
Developments in Mathematics, {\bf 30}. Springer, Cham, 2013.  

\bibitem[I]{I} N. It\^o, {\it Some Studies on Group Characters}. Nagoya Math. J.\,{\bf 2}, (1951), 17--28.

\bibitem[Ka]{Ka} I. Kaplansky, {\it Fields and Rings}. Chicago Lectures in Math., Univ. of Chicago Press, $1969$.

\bibitem[Ke]{K} O. H. Keller, 
{\it Ganze Gremona-Transformationen}. Monats. Math. Physik {\bf 47} (1939), no.\,1, 299-306. 

\bibitem[Ko]{Ko} A. Konijnenberg, {\it Mathieu Subspaces of Finite Products of Matrix Rings.} Master Thesis, Radboud University Nijmegen, 
The Netherlands, 2012.

\bibitem[Man]{Manz} O. Manz, {\it On the Modular Version of It\^o-Michler Theorem on Character Degrees for Groups of Odd Order}.  Nagoya Math. J. Vol.\,{\bf 105} (1987), 121--128. 


\bibitem[Mat]{Ma} O. Mathieu, {\it Some Conjectures about Invariant Theory and Their Applications.} Alg\`ebre non commutative, groupes quantiques et invariants (Reims, 1995), 263--279, S\'emin. Congr., 2, Soc. Math. France, Paris, 1997. 

\bibitem[Mi1]{M1} G.O. Michler, {\it Trace and Defect of a Block Idempotent}, J. Alg. {\bf 131}\,(1990) 496-501.

\bibitem[Mi2]{M2} G.O. Michler, {\it Brauer's Conjectures and the Classification of Finite Simple Groups}. Lect. Notes Math.\,{\bf 1178}, Springer, Berlin, 1986.

%
 
 

\bibitem[N1]{N} G. Navarro, {\it Characters and Blocks of Finite Groups}.  London Mathematical Society Lecture Notes Series, Cambridge University Press, 1998.
 
\bibitem[N2]{N2} G. Navarro, {\it 
Variations on the It\^o-Michler Theorem on Character Degrees}. 
Rocky Mountain J. Math. V.\,{\bf 46}, No. 4 (2016), 1363-1377.

\bibitem[Pa]{P} D.S. Passman, \newblock{\em The Algebraic Structure of Group Rings}, Pure and Applied mathematics, A Wiley-Interscience publication. John Wiley and Sons, Inc. 1977.


\bibitem[Pi]{Pi} R. S. Pierce, {\it Associative Algebras}. Graduate Texts in Mathematics, 88. Springer-Verlag, New York-Berlin, 1982. 


\bibitem[PD]{PD} B.M. Puttaswamaiah and  J.D. Dixon, \newblock{\em Modular representations of finite groups}. Pure and Applied mathematics, Academic Press Inc., 1977.

\bibitem[R]{R} G. R. Robinson, {\it On finite groups with all simple modules of low dimension in characteristic $p$}. 	arXiv:2009.10602 [math.GR].

\bibitem[Za]{Za} A. E. Zalesskii, {\em On a Problem of Kaplansky}. Soviet Math. {\bf 13} (1972), 449-452.
 

\bibitem[Z1]{IC} W. Zhao, {\it Images of Commuting  Differential Operators of Order One with Constant Leading Coefficients}.  J. Alg. {\bf 324} (2010),  no. 2, 231--247. [MR2651354]. See also arXiv:0902.0210 [math.CV]. 

\bibitem[Z2]{GIC} W. Zhao, {\em Generalizations of the Image Conjecture and the Mathieu Conjecture}.  J. Pure Appl. Alg. {\bf 214} (2010), 1200-1216. See also arXiv:0902.0212 [math.CV].  
 
\bibitem[Z3]{MS} W. Zhao, {\em Mathieu Subspaces of Associative Algebras}. J. Alg. {\bf 350} (2012), no.2, 245-272. 
See also arXiv:1005.4260 [math.RA].

\bibitem[ZW]{ZW} W. Zhao and R. Willems, \newblock {\em Analogue of the Duistermaat-Van Der Kallen Theorem for Group Algebras}. Cent. Eur. J. Math. {\bf 10}\,(3) (2012), 974-986.

\end{thebibliography}
\end{document}